
\input amstex
\documentstyle{amsppt}\magnification=\magstep1
\magnification=1200
\NoRunningHeads \NoBlackBoxes \baselineskip=12pt
\vsize=8.5truein \hsize=6.0truein \hoffset=0.3truein
\voffset=0.2truein \loadmsam \loadbold \loadmsbm \loadeufm
  \NoRunningHeads
\document
\define\({\left(}
\define\){\right)} \define\[{\left[} \define\]{\right]}
\define\db{\bar\partial} 
\define\aint{\rlap{\kern2pt\vrule height2.5pt width9.5pt depth-2.1pt}
\int_{B_{\rho}}} \TagsOnRight
     \nologo

\topmatter
\title  The $\db$-Cauchy problem   and nonexistence of Lipschitz  Levi-flat
hypersurfaces in $\Bbb CP^n$ with $n\ge 3$
         \endtitle
\author   Jianguo Cao and   Mei-Chi Shaw* \endauthor
\address
Department of Mathematics, University of Notre Dame, Notre Dame,
IN  46556 USA  \endaddress \email jcao$\@$nd.edu and
     Shaw.1\@nd.edu \endemail

\thanks{*Both authors are partially
supported by NSF grants.   }\endthanks

\endtopmatter

\centerline{October 12, 2006, To appear in Math. Zeit. }

\bigskip

\bigskip

In this paper we study the     $\db$-Cauchy problem   and the
$\db$-closed extension problem for forms  on  domains in complex
hermitian manifolds. These problems were first studied in the paper
by Kohn-Rossi \cite{KR} (see also \cite{FK}), who proved the
holomorphic
  extension of smooth $CR$  functions  and the $\db$-closed
extension of smooth forms from the boundary $b\Omega$ of a
strongly pseudoconvex domain to the whole domain $\Omega$. The
$L^2$ theory of these problem has been obtained for pseudoconvex
domains in $\Bbb C^n$ or, more generally, for domains in complex manifolds with
strongly plurisubharmonic weight functions (see Chapter 9 in
\cite{CS} and the references therein).  In this paper we study
these problems on pseudoconvex domains   in complex
hermitian manifolds when such weight functions   are not
available,  for instance, on a pseudoconvex domain in the complex projective
space
$\Bbb CP^n$.

One   application  of the $\db$-Cauchy problem is to obtain the
nonexistence of Levi-flat
hypersurfaces in $\Bbb CP^n$. This was first used by   Siu in \cite{Si1} where
the nonexistence of smooth (or $\frac {3n}2+7$) Levi-flat
hypersurfaces in $\Bbb CP^n$ was
proved for  $n\ge 3$. In a subsequent paper \cite{Si2}, he proved the
nonexistence
of $C^8$ Levi-flat hypersurfaces in $\Bbb CP^2$.
      We also mention the papers by Lins-Neto \cite{LN}, Iordan \cite{Io}
and Ni-Wolfson \cite{NW} on related subjects.

   The main result of this paper is to
   prove the nonexistence of Lipschitz Levi-flat hypersurfaces in
$\Bbb CP^n$ for $n\ge 3$.
      We first define Lipschitz Levi-flat hypersurfaces.

Recall that a  bounded domain ${D}\subset\subset\Bbb R^{2n}$ is called
Lipschitz if  near every boundary point   $p\in bD$, there
exists a neighborhood $U$
of $p$ such that in local coordinates
$(x',x_n)=(x_1,\cdots,x_{2n-1},x_{2n})$,   $${D}\cap U=\{(x',x_{
2n})\in
U\mid  x_{ 2n}>\psi(x')\} $$
for some Lipschitz function $\psi:\Bbb R^{2n-1} \to  \Bbb R $.   A Lipschitz
function is
differentiable almost everywhere
  (See Evans-Gariepy \cite{EG} for a proof of this fact). A domain in
a complex manifold
is called Lipschitz if at every point of the boundary, there exist
some local coordinates such that the
boundary is the graph of some Lipschitz function.

\proclaim{Definition}A Lipschitz hypersurface is  a hypersurface
which  locally is the graph of a
Lipschitz function.  A Lipschitz (or $C^1$)  hypersurface is  said to be
Levi-flat  if  it is locally
     foliated by complex manifolds of complex dimension $n-1$.
\endproclaim
 From the implicit function theorem, any $C^1$ hypersurface locally is
the graph of some $C^1$ function.
A $C^2$ hypersurface $M$ is called
Levi-flat if its Levi-form vanishes on $M$.
    Any $C^k$ Levi-flat
hypersurface, $k\ge 2$ is locally  foliated by complex manifolds
of complex dimension
$n-1$. The foliation is of class $C^{k}$ if the hypersurface is of
class $C^k$, $k\ge 2$ (see
    Barrett-Fornaess
\cite{BF}).   The proof in \cite{BF}
also gives that if a real $C^1$ hypersurface admits a continuous foliation
by complex manifolds, then the foliation
is actually $C^1$. Thus our definition is a natural generalization of
Levi-flatness to Lipschitz or $C^1$  hypersurfaces.

\proclaim{Theorem}    There exist  no Lipschitz  Levi-flat
hypersurfaces in $\Bbb
CP^n$ for $n\ge 3$.
\endproclaim

   The main tool to prove the theorem  is
to study the
$\db$-Cauchy problem using the $\db$-Neumann operator. When the
boundary is $C^2$ and pseudoconvex in $\Bbb
CP^n$,  the
$\db$-Neumann operator exists   using  bounded
plurisubharmonic functions, a result by Ohsawa-Sibony \cite{OS}.  It
is not known if the $\db$-Neumann operator exists for Lipschitz
pseudoconvex domains. However,
the weighted $\db$-Neumann operator always exists with suitable
weight functions.  To
  prove  the nonexistence of
Lipschitz Levi-flat hypersurfaces, we
    use  the  $L^2$ $\db$-Cauchy problem with weights and the
equivalence of the weighted spaces with the
Sobolev
    spaces.

In    \cite{CSW},  we carried out an $L^2$ approach for $\db$-closed
extension problem
  using the $\db$-Neumann operator  in order  to study the
nonexistence   of $C^2$-smooth Levi-flat
real hypersurfaces in $\Bbb CP^n$. In fact, only the
    nonexistence
of $C^{2,\alpha}$    Levi-flat hypersurfaces in $\Bbb CP^n$  with
$n\ge 3$  was proved,  by   using $\db$-closed extension of the Chern
connection
$(0,1)$-forms (see Section 5 in \cite{CSW}).   The proof for the
$\Bbb CP^2$ case in Section 6 of  \cite{CSW}  relies on a
Liouville-type result , which is yet
to be completed (see Conjecture 2 at the end of this paper).
At the end of the paper, we mention how to bridge the gap in the
proof \cite{CSW} for the nonexistence of $C^2$
Levi-flat hypersurfaces in
$\Bbb CP^2$ using results in \cite{Si2}.

  We note that there exist nonsmooth Levi-flat hypersurfaces in $\Bbb
CP^n$ which are not locally
Lipschitz graphs.  Let
$ M =\{[z_0,z_1,z_2]\in\Bbb CP^2\mid |z_0|=|z_1|\} $ and $\Bbb
CP^2\setminus M=\Omega^+\cup
\Omega^-$,
   where
$[z_0,z_1,z_2]$ are homogeneous coordinates in $\Bbb CP^2$.
  Then $\Omega^+ $ and $\Omega^-$ are both pseudoconcave and
pseudoconvex domains  since
each can be represented in local
coordinates   by a  product of a disc with $\Bbb C$ (see e.g.
\cite{HI}).   We can view $M$ as a
Levi-flat hypersurface in the sense that  it is the boundary of a
domain which is both    pseudoconvex and
pseudoconcave.   The boundary
$M$ is   smooth  except at $[0,0,1]$, where $M$ is   not foliated by
complex curves. Notice
that
$M$ is also not a graph of a Lipschitz function  in a neighborhood of
the point $[0,0,1]$.
   Similar examples can be found in $\Bbb CP^n$ for $n\ge 3$ by
setting $ M =\{[z_0,z_1,z_2, \cdots,
z_n ]\in\Bbb CP^n \mid |z_0|=|z_1|\} $.

The plan of this paper is as follows:  In section 1 we give a
self-contained treatment of  the
$\db$-Cauchy problem on   domains with Lipschitz boundary in a
hermitian complex
manifold  using   the
$\db$-Neumann operators.  In section 2 we prove the existence of
H\"older continuous bounded exhaustion
functions for pseudoconvex domains with $C^{1,1}$   boundary in $\Bbb CP^n$.
This gives an alternative proof of the Ohsawa-Sibony result
on the existence of bounded plurisubharmonic functions  for $C^2$ pseudoconvex
domains
in $\Bbb CP^n$.
In Section  3, we use the weighted
$\db$-Cauchy problem  to study the extension of
$\db$-closed
$(p,q)$-forms from a
pseudoconcave domain to $\Bbb CP^n$ when $q<n-1$, $n\ge 3$.    In Section 4, we
study the
Levi-flat boundary and its connection forms  and  prove the main theorem.
   It is still unknown if our main theorem can
be extended to
$\Bbb CP^2$. In Section 5, we discuss the extension of $\db$-closed
$(p,n-1)$-forms in $\Bbb CP^n$.
  We also mention  two open problems which will imply the
nonexistence of Lipschitz Levi-flat hypersurfaces in
$\Bbb CP^2$.

\heading
{\bf
1. The $L^2$ $\db$-Cauchy problem  on complex manifolds  }\endheading

Let $\Cal X$ be a   complex hermitian  manifold of dimension $n\ge 2$
and let $\Omega$ be
a bounded
domain in
$\Cal X$.
   The $L^2$  Cauchy problem for $\db$ is to study the following question:
   Given a $(p,q)$-form $f$ with $L^2$ coefficients   supported in
$\overline\Omega$, where $0\le p\le n$ and $1\le q\le n$,
    find a $(p,q-1)$-form $u$ such that
$$\cases &\text{Supp }u\subset\overline \Omega,
\\&\db u=f \quad \text{in } \Cal X  \text{ in the distribution
sense.}\endcases \tag 1.0  $$

When   $q<n$, we assume that $f$ satisfies
$$ \db f=0 \quad \text{in } \Cal X  \text{ in the distribution
sense}.\tag 1.1 $$
When $q=n$,    (1.1) is a void condition.  Using
integration-by-parts, another  compatibility
condition for
    (1.0) can be derived as follows: If (1.0)  is solvable for  $f\in
L^2_{(p,q)}(\Omega)$, where
$1\le q\le n$, then
$f$ must satisfy
$$
\int_\Omega f\wedge g=0,\quad g\in
L^2_{(n-p,n-q)}(\Omega)\cap\text{Ker}(\db).\tag 1.2
$$
We define the generalized Bergman projection operator
$$P_{(p,q)}:L^2_{(p,q)}(\Omega)\to
   L^2_{(p,q)}(\Omega)\cap \text{Ker}(\db).$$

Recall that the Hodge star operator $\star = \bar{*}$ is given by
   $$ (\star f,  g)|_\Omega = (-1)^{p+q}\overline{ \int_\Omega  g\wedge
f}= \overline{\int_\Omega f\wedge g}.  $$
Hence,   condition (1.2) is  equivalent to
$$P_{(n-p,n-q)}( {\star f)}=0. \tag 1.2'$$
Thus when $q<n$, both (1.1) and (1.2) are compatibility conditions for
the $\db$-Cauchy problem.

In the next lemma, we will show that
condition (1.2) implies condition (1.1).
\proclaim{Lemma 1.1} Let $\Omega$ be a bounded domain in a complex
hermitian manifold $\Cal X$ of dimension $n\ge 2$.
Let $f\in L^2_{(p,q)}(\Omega)$,  where $0\le p\le n$ and
    $0\le q<n$,    such that $f$ satisfies  (1.2). Then $\db f=0$ in $\Cal X$
if $f$ is extended to be zero outside $\Omega$.
\endproclaim
\demo{Proof}  We take $g=\db\star v$ for some $v\in
C^\infty_{(p , q+1)}(\Cal X)$ in (1.2). It is clear that $g \in
\text{Ker}(\db)$.
Let $ \vartheta  = \db^* = - \star \db \star$, where $\vartheta$ is
the formal adjoint of
$\db$ and $\db^*$ is the Hilbert space adjoint.   By (1.2) and the fact
$\db
\star v
\in
\text{Ker}(\db)$, we see that
$$(f, \db^* v)_{\Cal X}=\int_\Omega f\wedge \star (\db^*
v)=(-1)^{p+q+1}\int_\Omega f\wedge\db (\star v)=0$$ for any $v\in
C^\infty_{(p , q+1)}(\Cal X)$, where we used the equality $\star
(\star u) = (-1)^{p+q} u$ for $u = \db (\star v) \in
L^2_{(n-p,n-q)}(\Omega)$. This implies that $\db
f=0$ in the distribution sense in $\Cal X$.
\qed\enddemo

In general, (1.1) and (1.2) are not equivalent.
   We will see that they are equivalent for $q<n$
in Theorem 1.4.

When $q\le n$, including the top degree case,  the $\db$-Cauchy problem
will be solved for forms satisfying (1.2) in the next theorem.

\proclaim {Theorem 1.2}Let $\Omega$ be a bounded domain in a complex
hermitian manifold $\Cal X$ of dimension $n\ge 2$.  Suppose that  the
$\bar\partial$-Neumann operator $N_{(n-p,n-q)}$  on
$L^2_{(n-p,n-q)}(\Omega)$  exists  for some  $0\le p\le
n$ and
$1\le q\le n $.   For any
$f\in L^2_{(p,q)}(\Cal X)$
such that
$f$ is supported in $\overline \Omega$ and $f$ satisfies (1.2),
then   there exists  $u\in L^2_{(p,q-1)}(\Cal X)$ satisfies
$\bar\partial u=f$ in the
distribution sense in
$\Cal X$ with
$u$ supported in $\overline \Omega$.

\endproclaim
\demo{Proof}  Since the
$\bar\partial$-Neumann operators $N_{(n-p,n-q)}$  in
$\Omega$  exists,  the  generalized Bergman projection operator
   $P_{(n-p,n-q)}: L^2_{(n-p,n-q)}(\Omega)\to
L^2_{(n-p,n-q)}(\Omega)\cap \text{Ker}(\db)$  is given
by
$$\db^*\bar\partial N_{(n-p,n-q)}=I-P_{(n-p,n-q)}.\tag 1.3$$ We set $u$ by
   $$u=-\star {\bar\partial N_{(n-p,n-q)}{\star f}}.\tag 1.4 $$ Since
$f$ satisfies (1.2),
we have
$P_{(n-p,n-q)}\star f=0$.  From  (1.3),  we have
$$\aligned \db u&=(-1)^{p+q}
{\star\bar\partial^*\bar\partial N_{(n-p,n-q)}\star  f}\\&=
f-(-1)^{p+q} {\star P_{(n-p,n-q)}\star f}=f \quad \text{in
}\Omega.\endaligned \tag
1.5$$

Using the fact that $\star u\in \text{Dom}(\db^*)$ and extending $u$ to
be zero outside
$\Omega$,   one can show that   $\bar\partial u=f$ in $\Cal X$ in
the distribution
sense as follows.
    Observe that
   $$ \bar\partial^*(\star  u) =\vartheta\star  u  =(-1)^{p+q
}\star\bar\partial u =(-1)^{p+q }\star  f, $$ where
$\vartheta\star  u$ is taken in the distribution sense in
$\Omega$. Hence, we have   for any $\psi\in C^\infty_{(p,q)}(\Cal X),$
$$\aligned (u,\vartheta\psi)_{\Cal X} &=(\star
{\vartheta\psi},\ \star  u)_\Omega \\ &=(-1)^{p+q
}(\bar\partial\star \psi,\ \star  u)_\Omega \\ &=(-1)^{p+q
}(\star \psi,\ \bar\partial^*(\star  u))_\Omega \\ &=(\star
\psi,\ \star  f)_\Omega \\ &=(f,\ \psi)_{\Cal X},\endaligned \tag 1.6
$$ where the
third equality holds since $\star u\in
\text{Dom}(\bar\partial^*)$. Thus $\bar\partial u=f$ in the
distribution sense in $\Cal X$.
\qed
\enddemo

Theorem 1.2  implies that  condition (1.2) is necessary and  sufficient
for solving the
$\db$-Cauchy problem for all $(p,q)$-forms of all degrees, including
the top degree
$q=n$.

   Next we analyze the case when $q<n$.   Let
$\Cal H_{(p,q)}(\Omega)$ denote the space of harmonic
$(p,q)$-forms, i.e.,
$$\Cal H_{(p,q)}(\Omega)=\{h\in L^2_{(p,q)}{(\Omega)}\cap
\text{Dom}(\db)\cap \text{Dom}(\db^*)\mid
\db h=0,\
\db^*h=0\}.$$
Notice that no assumption on the smoothness of $\Omega$ is used in Lemma 1.1
and  Theorem  1.2. From now on, we will assume that the domain
$\Omega$ has  Lipschitz boundary.

\proclaim{Lemma 1.3}Let $\Cal X$ be a  complex hermitian
manifold of dimension
$n\ge 2$.  Let
$\Omega$ be a bounded  domain  in $\Cal X$ with Lipschitz
boundary.   For   $0\le p\le
n$,
   $1\le q\le n-1$, if $f\in L^2_{(p,q)}(\Cal
X)$  with
$\bar\partial f=0$ in the distribution sense in
$\Cal X$ and $f$ supported in $\overline\Omega$, then $\star f\in
\text{Dom}(\db^*)$ and
$\db^*\star f=0$ in $\Omega$.
\endproclaim
\demo{Proof}
For any
$\phi\in C^\infty_{(n-p,n-q-1)} (\overline\Omega) $,
$$ \aligned (\bar\partial\phi,\star  f)_\Omega  &
=(-1)^{p+q} \int_\Omega \db
\phi\wedge f=(-1)^{p+q}\int_\Omega f\wedge
\star\star \db\phi\\ &=(-1)^{p+q}(f,\star\db\phi)_\Omega =
(f, \vartheta\star \phi )_\Omega    \\ &= ( \bar\partial
f,\star \phi)_{\Cal X}\\ &=0 \endaligned $$
since $\text{supp}\  f\subset\overline \Omega$ and
$\bar\partial f=0$ in the distribution sense in $\Cal X$.

Since $\Omega$ has Lipschitz boundary $b\Omega$, using the
Friedrichs's lemma,   we see
that the set
$C^\infty_{(n-p,n-q-1)} (\overline\Omega) $ is dense in
Dom($\bar\partial$) in the
graph norm (see \cite{H\"o1} or  Step 1 in Lemma 4.3.2 in
\cite{CS}).   It follows
from the definition of
$\db^*$  that
$\star f\in\text{Dom}(\bar\partial^*)$ and
$\bar\partial^*(\star  f)=0$.\qed
\enddemo

We summarize the discussion above as follows.

\proclaim{Theorem 1.4}Let $\Cal X$ be a   complex hermitian
manifold of dimension
$n\ge 2$.  Let
$\Omega$ be a bounded   domain  in $\Cal X$ with Lipschitz
boundary.  We assume that  the
$\bar\partial$-Neumann operators $N_{(n-p,n-q)}$   and $N_{(n-p,n-q-1)}$   in
$\Omega$  exist  for   $0\le p\le n$ and
$1\le q\le n-1$ and  assume that $\Cal H_{(n-p,n-q)}(\Omega)=\{0\}$.
For every $f\in
L^2_{(p,q)}(\Cal X)$  with
$\bar\partial f=0$ in the distribution sense in
$\Cal X$ and $f$ supported in $\overline\Omega$, one can
find $u\in L^2_{(p,q-1)}(\Cal X)$ such that
$\bar\partial u=f$ in the distribution sense in $\Cal X$
with
$u$ supported in $\overline\Omega$.

\endproclaim
\demo{Proof} By our assumption,  the
$\bar\partial$-Neumann operator $N_{(n-p,n-q)}$ of degree  $(n-p,n-q)$ in
$\Omega$  exists and  $\Cal H_{(n-p,n-q)}(\Omega)=\{0\}$.  From the
Hodge decomposition,  we have
for every $f\in
L^2_{(p,q)}(\Omega)$,
$$\star f=\db\db^*N_{(n-p,n-q)}\star f+\db^*\db N_{(n-p,n-q)}\star
f.$$ We define
$$u=- {\star\bar\partial N_{(n-p,n-q)} {\star f}},\tag
1.7$$ then $u\in L^2_{(n-p,q-1)} (\Omega) $ and
$\star   u\in \text{Dom} (\bar\partial^*)$.

Extending $u$ to $\Cal X$ by defining $u=0$ in $\Cal X
\setminus \Omega$, we claim  that $\bar\partial u=f$ in the
distribution sense in $\Cal X$. First we  prove  that
$\bar\partial u=f$ in the distribution sense in $\Omega$.

  By (1.7) we get
$$ \aligned
\bar\partial u&=- {\bar\partial\star\bar\partial
N_{(n-p,n-q)}\star  f}\\ &=(-1)^{p+q+1}
{\star\star\bar\partial\star\bar\partial
N_{(n-p,n-q)}\star  f}\\ &=(-1)^{p+q}
{\star\vartheta\bar\partial N_{(n-p,n-q)}\star  f}\\
&=(-1)^{p+q} {\star\bar\partial^*\bar\partial
N_{(n-p,n-q)}\star  f}. \endaligned\tag 1.8$$

It follows from Lemma 1.3 that $\star f$ is in Dom($\db^*)$ and
$$\db^*(\star f)=0.\tag 1.9$$
By our assumption that $N_{(n-p,n-q-1)}$ exists, we have
$$\bar\partial^*N_{(n-p,n-q)}\star
f=N_{(n-p,n-q-1)}\bar\partial^*(\star  f)=0.\tag 1.10$$
Combining (1.8) and (1.10) and  the assumption $\Cal
H_{(n-p,n-q)}(\Omega)=\{0\}$, we conclude that
$$\aligned \bar\partial u&=(-1)^{p+q}\star
{\bar\partial^*\bar\partial N_{(n-p,n-q)}\star  f}\\
&=(-1)^{p+q}\star {(\bar\partial^*\bar\partial
+\bar\partial\bar\partial^*)N_{(n-p,n-q)}\star  f}\\
&=(-1)^{p+q}\star\star f\\ &=f \endaligned  $$ in the
distribution sense in $\Omega$.  Since
$\star u\in\text{Dom}(\bar\partial^*)$, repeating the same arguments
as in (1.6), we have
proved $\bar\partial u=f$ in the distribution sense in $\Cal X$.
Theorem 1.4 is proved.\qed
\enddemo

    We note that in the proof  of Lemma 1.3 and Theorem 1.4, the
Lipschitz boundary condition
on
$\Omega$ is used to show that the
$C^\infty_{(n-p,n-q-1)}(\overline\Omega)$ space is dense in
Dom($\db$) in the graph norm.

Let $\Omega\subset\subset\Bbb CP^n$
be a pseudoconvex domain with  $C^2$-smooth boundary
$b\Omega$ and let $\delta(x) = d(x, b\Omega)$ be the distance
function from $x\in\Omega$ to
$b\Omega$. We call $t_0 = t_0(\Omega) $   the
order of plurisubharmonicity
for
the distance function  $\delta$ if
$$ t_0(\Omega) = \sup \{ 0<\epsilon \le 1 | i\partial\bar{\partial} (-
\delta^{\epsilon   }) \ge 0 \text{ on } \Omega\}.  \tag 1.11 $$
     In
$\Bbb CP^n$ with the standard Fubini-Study metric,     Ohsawa-Sibony
[OS] showed  that there
exists $ 0< t_0(\Omega)\le 1  $ for any
pseudoconvex domain $\Omega \subset \Bbb
CP^n$ with $C^2$-smooth  boundary (see  Diederich-Fornaess \cite{DF}
for domains in
$\Bbb C^n$).  We recall the following results (see Theorem 2 in
     \cite{CSW}).

\proclaim{Theorem 1.5} Let $\Omega$ be a pseudoconvex domain   with
$C^2$-smooth boundary in
$ \Bbb CP^n$ and  let  $ t_0$ be the order of plurisubharmonicity
for the distance function $\delta$.
Then the $\db$-Neumann operator
$N_{(p, q)}$    exists on
$L^2_{(p, q) }(\Omega)$ where $0\le p,  q\le n$ and  the harmonic forms
$\Cal H_{(p,q)}(\Omega)=\{0\}$ if $1\le q\le n$.    Furthermore,  $N,
\bar{\partial}N,
\bar{\partial}^*N$ and the Bergman projection $   P$ are  exact regular
      on
$W^{s}_{(p, q) }(\Omega)$ for $0 \le s <
\frac 12 t_0 $ with respect to the  $W^s(\Omega)$-Sobolev norms.
\endproclaim

A direct consequence of  Theorems 1.2, 1.4 and 1.5 for the case of
$\Cal X=\Bbb CP^n$ is the corollary below,  which was
already  obtained  in   Propositions 4.1 and 4.2 in
\cite{CSW}.
\proclaim{Corollary 1.6 ($L^2$ Cauchy problem for $\db$ in $\Bbb CP^n$)}   Let
$\Omega\subset\subset \Bbb CP^n$ be a  pseudoconvex domain    with
$C^2$ boundary
and let
$0\le p\le n$ and
$1\le q\le n$.   For every
$f\in L^2_{(p,q)}(\Bbb CP^n)$   supported in $\overline\Omega$, we assume that
$\bar\partial f=0$ in the distribution sense in
$\Bbb CP^n$ if $1\le q\le n-1$ and $f$ satisfies (1.2) if $q=n$. Then one can
find $u\in L^2_{(p,q-1)}(\Bbb CP^n)$ such that
$\bar\partial u=f$ in the distribution sense in $\Bbb CP^n$
with
$u$ supported in $\overline\Omega$.

Furthermore, if  $f\in W^s_{(p,q)}(\Omega )$ with $0\le s<\frac 12
t_0$, then we can choose $u\in
W^s_{(p,q)}(\Omega )$.

\endproclaim

In the next section, we will show that when the domain is
pseudoconvex with $C^{1,1}$ boundary,  then
Theorem 1.5 and Corollary 1.6 hold.

\heading {\bf 2.  Bounded plurisubharmonic functions for pseudo-convex domains
with $C^{1, 1}$    boundary}\endheading

In this section we will  recall some results for pseudoconvex domains in $\Bbb
CP^n$. We
will also give an alternative proof
of the existence of bounded plurisubharmonic functions for domains with $C^{1,
1}$
boundary (see \cite{OS}). Such functions can be used to prove the existence of
the $L^2$
$\db$-Neumann operators.

\proclaim{Lemma 2.1} Let $\Omega$ be a Lipschitz pseudoconvex domain with
Levi-flat boundary $M$ in $\Bbb CP^n$, $n\ge 2$. Then $M$ is locally
foliated by complex hypersurfaces. Moreover,
 for each $Q\in M$,
there exist a neighborhood $U $ of $Q$ and local unitary frame $\{\tilde{e}_1,
...., \tilde{e}_{n-1}, \tilde{e}_n\}$
 on  $U $    such that (1) for $z \in M \cap U$,
 the vector fields  $  \{\tilde{e}_1, ...., \tilde{e}_{n-1}\}|_z $ are tangent
to the leaves
of the foliation of $M\cap U$; and (2) The covariant derivative $\nabla_\xi
\tilde{e_j}$  is a
bounded function for $j=1, ...,  n-1 $ and any unit vector
$\xi \in T_z(\Bbb CP^n  )$ with $z \in U$.

\endproclaim
\demo{Proof}
 Since $M$ is Levi-flat, it is locally foliated
by complex manifolds of dimension $n-1$ and  the foliation is Lipschitz in  the
transversal direction.
 For any point
$Q \in M$, we can parametrize a neighborhood $V\subset M$ of $Q$ as follows.
Let $\{z' ,g(z' ,t)\}$ denote the leaf $\Sigma_t$ where $g(z'
,t)$ is holomorphic in $z'=(z_1,\cdots,z_{n-1}) \in \Bbb B_{\epsilon} \subset
\Bbb C^{n-1} $  and Lipschitz  in $t$ for   $0\le |t|<\mu$. We can
parametrize $M$ locally as a graph of the function $ g$, by setting
$$\Psi(z' ,t)=(z' , g(z' ,t)),
$$ where
$z' \in \Bbb C $,
     $0\le |t|<\mu$.  Clearly, $\Psi: \Bbb B_{\epsilon} \times (-\mu, \mu)
\to M$ is a local coordinate map of $M$ and $\Psi$ is Lipschitz in $t$ and
$C^\infty$ (holomorphic) in $z' $.

  Let $z' = (z_1, ..., z_{n-1})$ and extend $\Psi$ to a map $\tilde \Psi:
   \Bbb B_{\epsilon} \times (-\mu, \mu) \times (-\mu, \mu)
\to  \Bbb CP^n $ by setting $ \tilde \Psi(z', t + is) = (z' , g(z'
,t) + s \vec{v}_0), $
 where $(t, s) \in (-\mu, \mu) \times (-\mu, \mu)$ and $\vec{v}_0$ is a constant
vector transversal to $ \frac{\partial g(z'
,t)}{ \partial t}  $ for all $(z',  t) \in \Bbb B_{\epsilon} \times (-\mu,
\mu)$. We now  choose
$\tilde{v}_j = \frac{\partial \tilde \Psi}{\partial z_j}$ for $j = 1, ...,
 n-1 $. Applying the Gram-Schmidt process to the frame $\{\tilde v_1, ...,
\tilde{v}_{n-1}\}$, we   obtain a unitary frame $\{ \tilde{e}_1, ....,
\tilde{e}_{n-1}   \}$ with Lipschitz coefficients. Thus (2) is satisfied  as
desired.
\qed

\enddemo
We recall the following theorem by \cite{Ta} (see also  \cite{CS}).

\proclaim{Theorem 2.2} Let $\Omega\subset\subset\Bbb CP^n$
be a pseudoconvex domain.   Then the distance function $\delta$    satisfies
  $$i\partial\db (-\log \delta )\ge  \omega   \tag 2.1$$
  as currents where $\omega$ is the K\"ahler form of the Fubini-Study
metric on $\Bbb CP^n$.

\endproclaim

In $\Bbb CP^n$ with the standard Fubini-Study metric,     Ohsawa-Sibony [OS]
showed  that there exists a bounded plurisubharmonic functions for
pseudoconvex domains with $C^2$ boundary. We give a proof below for
pseudoconvex domains with
$C^{1, 1}$ boundary.

\proclaim{Proposition 2.3  } Let $\Omega\subset\subset\Bbb CP^n$ be a
pseudoconvex domain with $C^{1, 1}$  boundary $b\Omega$. Then there exists
a distance function $\delta$   in $C^{1,1}(\overline\Omega)$ which satisfies
(2.1)
almost everywhere.
   Furthermore, there exists  $t_0 = t_0(\Omega)  $ with $0<t_0\le 1$
such that
$$
i\partial\db(-  \delta^{t_0})\ge 0. \tag 2.2$$

\endproclaim
 \demo{Proof}
Let $\delta$ be the distance function from $z\in  \Omega$ to $b\Omega$. Since
the boundary is of class $C^{1,1}$, we have that  there exists a
neighborhood $U$ of $b\Omega$ such that $\delta$ is in
$C^{1,1}(\overline\Omega\cap
U)$.
Using \cite{Ta}, we have
$$  i\partial\db(-\log \delta)=i\frac{\partial\db(- \delta)}{ \delta}+ \frac
{i\partial \delta\wedge\db \delta}{ \delta^2}
   \ge \omega \tag 2.3$$
near the boundary almost everywhere.

To prove (2.2),  observe that inequality (2.2) is equivalent to
$$i\frac{\partial\db(- \delta)}{ \delta}+(1-t_0)\frac
{i\partial \delta\wedge\db \delta}{ \delta^2}\ge 0.\tag 2.4$$ Compare (2.4) with
(2.3),
  we see that  (2.2) is equivalent to
$$
i\partial\db(-\log \delta)\ge t_0 \frac {i\partial  \delta\wedge\db \delta}{
\delta^2}.\tag 2.5
$$
Near a boundary point, we choose a special orthonormal basis $w_1,\cdots,w_n$
for $(1,0)$-forms such that $w_n=\sqrt 2\partial(-\delta)$. Let
$L_1,\cdots,L_n$ be its dual and let $a$ be any $(1,0)$-vector.  We decompose
$a=a_\tau+a_\nu$ where $a_\nu=\langle a,L_n\rangle$ is the complex
normal component and $a_\tau$ is the complex tangential component.  We have
$$
\aligned  &\langle  \partial\db(-\log \delta),a\wedge\bar a\rangle \\&=\langle
\frac{ \partial\db(-  \delta)}{ \delta},a_\tau\wedge\bar
a_\tau\rangle +2\Re\langle \frac{ \partial\db(-  \delta)}{
\delta},a_\tau\wedge\bar a_\nu\rangle
\\&+\langle
\frac{ \partial\db(-  \delta)}{ \delta},a_\nu\wedge\bar a_\nu\rangle +\frac {
|a_\nu|^2}{ \delta^2}.
\endaligned\tag 2.6
$$
From (2.1) and (2.3),  we have
$$\langle
 \partial\db(- \log  \delta)  ,a_\tau\wedge\bar a_\tau\rangle\ge \langle
\frac{ \partial\db(-  \delta)}{ \delta},a_\tau\wedge\bar a_\tau\rangle\ge
|a_\tau|^2  . $$ Thus from (2.6),
$$
 \aligned\langle  \partial\db(-\log \delta),a\wedge\bar a\rangle & \ge
 |a_\tau|^2 + \frac { |a_\nu|^2}{ \delta^2} -2| \langle
\frac{ \partial\db(-  \delta)}{ \delta},a_\tau\wedge\bar a_\nu\rangle |
\\& -|\langle
\frac{ \partial\db(-  \delta)}{ \delta},a_\nu\wedge\bar a_\nu\rangle |.
 \endaligned\tag 2.7
$$
  Using the
assumption that $b\Omega$ is $C^{1,1}$,  we have

$$ |\partial\db \rho|\le C \tag 2.8$$   Also for any $\epsilon>0$,  there exists
a small
neighborhood
$U$ of $b\Omega$   such that
$$
  |\partial\db \delta  |\le
 \frac{\epsilon}{\delta}.
 \tag 2.9
$$

   Thus for any $\epsilon>0$, we have from
(2.8),
$$|\langle
\frac{ \partial\db(-  \delta)}{ \delta},a_\tau\wedge\bar a_\nu\rangle |\le
C\(\frac 1\epsilon |a_\tau|^2+\epsilon \frac{|a_\nu|^2}{ \delta^2}\),
\tag 2.10$$
 and from (2.9),
$$|\langle
\frac{ \partial\db(-  \delta)}{ \delta},a_\nu\wedge\bar a_\nu\rangle| \le  \frac
{  \epsilon}{ \delta^2} |a_\nu|^2 \tag 2.11$$ on  a
sufficiently small neighborhood $U$ of the boundary.

Substituting (2.9)-(2.10) into (2.7) and choosing $\epsilon$ sufficiently small,
we have
$$\langle  \partial\db(-\log \delta),a\wedge\bar a\rangle\ge \frac 12 \frac
{|a_\nu|^2}{ \delta^2}
-K|a_\tau|^2\tag 2.12$$ for some large constant $K$ depending on $\epsilon$.
Multiplying (2.1) by $K$ and adding it to (2.12), we have
$$(K+1)\langle  \partial\db(-\log \delta),a\wedge\bar a\rangle\ge \frac 12 \frac
{|a_\nu|^2}{ \delta^2}. $$ This proves (2.5) with $t_0=\frac 1{2(K+1)}$ near the
boundary, or equivalently,  (2.2) is proved near the boundary.
Since $\Omega$ is Stein, on any relatively compact submanifold
$\Omega'\subset\subset \Omega$,  there exists a bounded strictly
plurisubharmonic
function on $\overline\Omega'$. By standard arguments one can extend $\delta$ so
that $\delta$ is the distance function near the boundary and
$\delta$ satisfies (2.1) and (2.2) in $\Omega$. \qed
\enddemo
\noindent
{\bf Remark:} Diederich-Fornaess \cite{DF} show that if $\Omega$ is a
pseudoconvex domain in $\Bbb C^n$ with $C^2$ boundary,   let
$\tilde\delta=\delta e^{-K|z|^2}$ with large $K>0$ , then (2.1) holds with
$\delta$
substituted by $\tilde\delta$.  The proof of Proposition 2.3 is a modified proof
of the Diederich-Fornaess \cite{DF} and Ohsawa-Sibony \cite{OS}
results. We also remark that bounded plurisubharmonic exhaustion functions exist
for pseudoconvex domains  in $\Bbb C^n$ with $C^1$ (see
Kerzman-Rosay \cite{KeR}) or even Lipschitz boundary
  (see Demailly \cite{De}), but it is not known if    such functions exist for
$C^1$ or Lipschitz pseudoconvex domains in $\Bbb CP^n$.

\proclaim{Proposition 2.4}
Let $\Omega$ be a  pseudo-convex domain   with
$C^{1,1}$-smooth  boundary in
$ \Bbb CP^n$, $n\ge 2$.
Then the $\db$-Neumann operator
$N_{(p, q)}$    exists on
$L^2_{(p, q) }(\Omega)$ where $0\le p,  q\le n$ and  the harmonic forms
$\Cal H_{(p,q)}(\Omega)=\{0\}$ if $1\le q\le n$.    Furthermore,
there exist $t_0>0$ such that  $N,
\bar{\partial}N,
\bar{\partial}^*N$ and the Bergman projection $   P$ are  exact regular
      on
$W^{s}_{(p, q) }(\Omega)$ for $0 < s <
\frac 12 t_0   $ with respect to the  $W^s(\Omega)$-Sobolev norms.
\endproclaim

\demo{Proof}
Let $\Omega$, $\delta$ and $t_0$  be the same as in Proposition 2.4.
  The proposition  follows exactly the same as the proof of Theorem 2 in
\cite{CSW}.

\enddemo

 From Proposition 2.4,
   the results of Theorem 1.5 and Corollary 1.6 hold also for
$C^{1, 1}$ pseudoconvex domains.
Then  we can use the same arguments as in Section 5 in  \cite{CSW} to show the
nonexistence of
$C^{1, 1}$ Levi-flat hypersurfaces in $\Bbb CP^n$ when $n\ge 3$.
But it is not known if  Proposition 2.3   holds for Lipschitz
domains.  In the next section,  we
will use
  the weighted $\db$-Neumann operators to study the $\db$-Cauchy
problem on Lipschitz domains.

\heading
{\bf
3.  The $\db$-Cauchy problem with weights on Lipschitz pseudoconvex
domains   in
$\Bbb CP^n$
}\endheading

  Let $\Omega $ be a pseudoconvex domain with Lipschitz boundary  in
$\Bbb CP^n$, $n\ge 2$.
We study the $\db$-Cauchy problem with weights and the $\db$-closed
extension of forms  from
   pseudoconcave domains.

For $t>0$, let
$L^2(e^{-\phi_t},\Omega)=L^2(\delta^t,\Omega)=L^2(\delta^t)$   be the
weighted $L^2$ space with
respect to the weight function
$\phi_t=-t\log
\delta$. The norm in $L^2(\delta^t)$ is denoted by $\|\ \|_{(t)}$.
Let $\db$ and $\db^*_t$ be the
closure of
$\db$ and its
$L^2$ adjoint with respect to the weighted $L^2(\delta^t)$ space.

\proclaim{\bf Proposition 3.1} Let $\Omega\subset\subset\Bbb CP^n$
be a pseudoconvex domain.
For any $t>0$ and $(p,q)$-form
   $f\in L^{2} (\delta^t )
$, where
$0\le p\le n$ and
$1\le q\le n$, such that
$\db f=0$ in
$\Omega$,
    there exists $u\in L^2_{(p,q-1)}(\delta^t)$ satisfying $\db u=f$ and
$$\|u\|_{(t)}^2 \le  \frac 1{t}\|f\|^2_{(t)}. \tag 3.1  $$
Furthermore, the weighted  $\db$-Neumann operator $N_t$ exists for all $t>0$.

\endproclaim
\demo{Proof} We first assume that $\Omega$ is $C^2$.  By  \cite{Ta},
we have that $\phi=-\log \delta$
is strictly plurisubharmonic
and $i\partial\db\phi\ge   \omega$, where $\omega$ is the K\"ahler
form of $\Bbb CP^n$ with the Fubini-Study metric.
   Using  H\"ormander's weighted
$L^2$ estimates for the
$\db$-Neumann problem (see e.g. Proposition A.4 in \cite{CSW}),
  we have  the
following   formula: for any $(p,q)$-form $g\in
\text{Dom}(\db)\cap\text{Dom}(\db^*_{t })$,
$$
\|\db g\|^2_{(t)}+\|\db^*_{t
} g\|^2_{(t)}\ge t((i\partial
\db  \phi) g,\bar g)_{(t)}. \tag 3.2
$$

  Thus, we have
$$\|\db g\|^2_{(t)}+\|\db^*_{t } g\|^2_{(t)}\ge t\|g\|^2_{(t)}. \tag 3.3$$

For any $f  \in L^{2}(\delta^t )$, there exists $u\in L^2(\delta^t)$
satisfying $\db u=f$ and
  (3.1). This proves the proposition when $\Omega$ is $C^2$. The
general case follows from
approximating the domain $\Omega$ from inside by smooth pseudoconvex
domains.\qed

  From (3.1), we have that the weighted $\db$-Neumann operator $N_t$
exists for each $t>0$ (see the proof of
Theorem 4.4.1 in \cite{CS}). \qed

\enddemo

We remark that there is  no smoothness assumption on the boundary
$b\Omega$ in Proposition 3.1. We
will use the weighted $\db$-Neumann operator $N_t$ to study the
$\db$-Cauchy problem.

\proclaim {Proposition 3.2}Let $\Omega\subset\subset\Bbb CP^n$
be a pseudoconvex domain with Lipschitz boundary, $n\ge 3$.
Suppose that
$f\in L^2_{(p,q)}(\delta^{-t},\Omega)$   for some $t>0$,  where   $0\le p\le
n$ and
$1\le q< n $.
Assuming that
$\db f=0$  in $\Bbb CP^n$ with $f=0$ outside $\Omega$,
then   there exists  $u_t\in L^2_{(p,q-1)}(\delta^{-t},\Omega)$  with
$u_t=0$ outside $\Omega$ satisfying
$\bar\partial u_t=f$ in the
distribution sense in
$\Bbb CP^n$.

\endproclaim
\demo{Proof}  From Proposition 3.1,  the weighted
$\bar\partial$-Neumann operators $N_t$
    exists for forms in $L^2_{(n-p,n-q)}(\delta^t,\Omega)$.
Let $\star_{(t)}$ denote the Hodge-star operator with respect to the
weighted norm
$L^2(\delta^t,\Omega)$. Then $$\star_{(t)}=\delta^t\star=\star
\delta^t$$ where $\star$ is the Hodge star
operator with the  unweighted $L^2$ norm.
Since
$f\in  L^2_{(p,q)}(\delta^{-t},\Omega)$, we have that
  $\star_{(-t)}f\in   L^2_{(p,q)}(\delta^{ t},\Omega).$
   Let  $u_t$ be defined  by
   $$u_t=-\star_{(t)} {\bar\partial N_t{\star_{(-t)} f}}.\tag 3.4 $$
Then $u_t \in   L^2_{(p,q-1)}(\delta^{-t},\Omega)$, since
${\bar\partial N_t{\star_{(-t)} f}}$ is in
$\text{Dom}(\db^*_t)\subset  L^2_{(n-p,n-q+1)}(\delta^{ t},\Omega)$.
Since $\db^*_t=\delta^{-t}\vartheta\delta^t=-\star_{(-t)}\db\star_{(t)}$,
using the same proof as in Lemma 1.3, we have  $\star_{(-t)}  f\in
\text{Dom}(\db^*_t)$ and $\db^*_t
\star_{(-t)}  f=0$ in
$\Omega$.
This gives
$$\bar\partial^*_t N_t\star_{(-t)}  f=N_t\bar\partial^*_t
\star_{(-t)}  f=0.\tag 3.5$$
 From (3.5), we have
$$\aligned \db u_t&=-\db( \star_{(t)} {\bar\partial N_t{\star_{(-t)}
f}})\\&=(-1)^{p+q}
{\star_{(t)}\bar\partial^*_t\bar\partial N_t\star_{(-t)}  f}
\\&=(-1)^{p+q}
{\star_{(t)}\bar\partial^*_t\bar\partial N_t\star_{(-t)}  f}+ (-1)^{p+q}
{\star_{(t)}\bar\partial\bar\partial^*_t  N_t\star_{(-t)}  f}\\& =(-1)^{p+q}
{\star_{(t)} \star_{(-t)}  f}
\\&=
f \quad \text{in
}\Omega.\endaligned \tag
3.6$$

First notice that   $\star_{(-t)}(-1)^{p+q}\db N_t\star_{(-t)}f=\db
N_t\star_{(-t)}f\in
\text{Dom}(\db^*_t)$.
We also have
$ \db^*_t \star_{(-t)} u=(-1)^{ p+q} \star_{(-t)}f$ in $\Omega$.
Extending $u_t$ to
be zero outside
$\Omega$, one can show that
$\bar\partial u_t=f$ in
$\Bbb CP^n$. The proof is similar to  the proof of Theorem 1.2.
    In fact, for any $\psi\in C^\infty_{(p,q)}(\Bbb CP^n),$
$$\aligned (u,\vartheta\psi)_{\Bbb CP^n} &=(\star
{\vartheta\psi},\ \star_{(-t)}  u)_{( t)\Omega} \\ &=(-1)^{p+q
}(\bar\partial\star \psi,\ \star_{(-t)}  u)_{( t)\Omega} \\ &=(-1)^{p+q
}(\star \psi,\ \bar\partial^*_t(\star_{(-t)}  u))_{( t)\Omega} \\ &=(\star
\psi,\ \star_{(-t)}  f)_{(t)\Omega}= (\star
\psi,\ \star   f)_{ \Omega}\\ &=(f,\ \psi)_{\Bbb CP^n},\endaligned \tag 3.7
$$ where the
third equality holds since $\star_{(-t)} u\in
\text{Dom}(\bar\partial^*_t)$. Thus $\bar\partial u=f$ in the
distribution sense in $\Bbb CP^n$.
\qed
\enddemo

\proclaim{Theorem 3.3} Let $\Omega \subset\subset \Bbb CP^n$ be a
pseudoconvex domain with Lipschitz
boundary and let $\Omega^+=\Bbb CP^n\setminus \overline\Omega$.
   For any $f\in W^{1+\epsilon}_{(p,q)}({\Omega}^+)$, where $0\le p\le
n$, $0\le q< n-1$ and
  $0<\epsilon<\frac 12$,
such that  $\db f = 0$ in $\Omega^+$, there exists
$F\in W^\epsilon_{(p,q)}(\Bbb CP^n)$ with
$F|_{\Omega^+}=f$ and  $\db F=0$ in $\Bbb CP^n$ in the
distribution sense.

\endproclaim
\demo{Proof}
   Since $\Omega$ has Lipschitz  boundary, there exists a bounded extension
operator from
   $W^s(\Omega^+)$ to $W^s(\Bbb CP^n)$ for all $s\ge 0$ (see e.g.
\cite{Gr} or \cite{St}).  Let
$\tilde f\in W^{1+\epsilon}_{(p,q)}(\Bbb CP^n)$  be the  extension of
$f$  so that
$ \tilde f|_{ {\Omega}^+  } = f $
with
   $\|\tilde f\|_{W^{1+\epsilon}(\Bbb CP^n)}\le
C\|f\|_{W^{1+\epsilon}(\Omega^+)}.$
Furthermore, we can choose an extension such that $\db \tilde f\in
W^\epsilon (\Omega)\cap
L^2(\delta^{-2\epsilon} ,\Omega)$.

   We define $T\tilde f$   by
   $T\tilde f = - {\star_{(2\epsilon)}\bar\partial N_{2\epsilon }
{(\star_{(-2\epsilon)}  \db \tilde
f)}}$ in
$ \Omega$.
      From
Proposition 3.2, we have that $T\tilde f\in
L^2(\delta^{-2\epsilon},\Omega)$. But for a Lipschitz domain, we
have that $T\tilde f\in L^2(\delta^{-2\epsilon},\Omega)$ is
comparable to $W^{\epsilon}(\Omega)$ when
$0<\epsilon<\frac 12$.
This gives that
$T\tilde f\in W^{\epsilon}(\Omega)$ and $T\tilde f$
   satisfies
   $\db
T\tilde f=\db \tilde f\quad \text {in }  \Bbb CP^n $
in the distribution sense  if we extend $T\tilde f$ to be zero
outside $\Omega$.

Since $0<\epsilon<\frac 12$, the extension by 0 outside $\Omega$  is
a  continuous operator from
$W^\epsilon(\Omega)$ to
$W^s(\Bbb CP^n)$ (see e.g. \cite{LM} or \cite{Gr}).   Thus we have
$T\tilde f\in
W^\epsilon(\Bbb CP^n)$.

Define
$$F=\cases &f,\quad x \in
\overline{\Omega}^+,\\&\tilde f-T\tilde f ,\quad x\in
\Omega.\endcases $$
Then $F\in W^\epsilon_{(p,q)}(\Bbb CP^n)$ and $F$ is a $\db$-closed
extension of $f$.
\qed

\enddemo

\proclaim{Corollary 3.4}
  Let $\Omega^+$ be a pseudoconcave domain
in $\Bbb CP^n$ with Lipschitz boundary, where $n\ge 2$.   Then
$W^{1+\epsilon}_{(p,0)}( \Omega^+)\cap \text{Ker}(\db)=\{0\}$  for
every $1\le  p\le n$
and  $W^{1+\epsilon}_{(0,0)}( \Omega^+)\cap \text{Ker}(\db)=\Bbb C$.
\endproclaim
  \demo{Proof}  Using Theorem 3.3 for $q=0$, we have   that any
holomorphic $(p,0)$-form on $\Omega^+$
extends to be a holomorphic
$(p,0)$ in $\Bbb CP^n$, which are zero (when $p>0$) or constants (when $p=0$).
\enddemo
\proclaim{Corollary 3.5} Let $\Omega^+$ be a pseudoconcave domain
in $\Bbb CP^n$ with Lipschitz boundary, where $n\ge 3$.
   For any $f\in W^{1+\epsilon}_{(p,q)}({\Omega}^+)$, where $0\le p\le
n$, $1\le q< n-1$, $p\neq q$ and
   $0<\epsilon<\frac 12$,
such that  $\db f = 0$ in $\Omega^+$, there exists
$u\in W^{1+\epsilon}_{(p,q-1)}(\Omega^+)$ with
$\db u=f$ in $\Omega^+$.

\endproclaim

\heading {\bf 4.   Nonexistence
of Lipschitz  Levi-flat hypersurfaces in $\Bbb CP^n$ when $n\ge 3$}\endheading

     In this section we study $\db_b$-exactness of  (0, 1)-form $f$ on
a Lipschitz
Levi-flat hypersurface $M \subset \Bbb CP^n$ and prove the
main theorem.  It is a  refinement  of arguments used in \cite{Si1}
and \cite{CSW}.

   We recall the definition of
the Chern connection form for the complex line bundle generated by
the complex normal of $M$. Let $\Bbb CP^n\setminus M=\Omega^+\cup\Omega^-$.
     Let $\rho$ be the signed distance function of $M$
$$\rho (z) =\cases  - d(z, M) ,
\qquad &\text{if}\   z\in \Omega^-, \\  d(z, M),\qquad \
\ \ &\text{if}\   z\in \overline{\Omega^+}.\endcases\tag 4.1$$

   If $J$ is the complex structure of $\Bbb CP^n$ and $\nabla$ is the covariant
derivative of $\Bbb CP^n$
with respect to the Fubini-Study metric, the   connection form of the
complex normal line bundle $\nabla \rho \otimes \Bbb C$ on
$M$ is given by
$$
\beta ( X) = \langle \nabla_{X} (\nabla \rho), J\nabla \rho \rangle = -\langle
\nabla_{X}(J\nabla \rho), \nabla \rho \rangle, \tag4.2
$$
where $X$ is a tangent vector on $M$
(see   (5.3) and (A.7) in \cite{CSW}).

  For a general hypersurface, we need  $C^2$ smoothness to define the
curvature form and the connection
form. In this case,  the curvature form $\tilde\Theta^N$ associated with
the complex line bundle for $M$
  is  a well-defined 2-form with   $C^0$ coefficients in
$U$ and    is
$d$-exact.
  Following the Chern formula (see Proposition A.1 in \cite{CSW}), we have that
$\tilde\Theta^N= \sqrt{-1}d\beta$  on a
tubular neighborhood  $U (M)$
of
$M$ in
$\Bbb CP^n$.
  Let  $\beta_b$ be the projection of $\beta$ to $M$ defined by
  $$\beta_b = \beta |_{T^{(1, 0)}(M)\oplus T^{(0, 1)}(M) }.$$
    Write
   $
\beta_b=\beta_b^{1,0}+\beta_b^{0,1}$  where
$\beta_b^{1,0}$ and
$\beta_b^{0,1}$ are the (1,0) and (0,1) components of $\beta_b$.
When the hypersurface $M$ is Levi-flat, one can relax the smoothness
using Lemma 2.1.  We first show that the
Chern connection and the curvature can be defined for
Lipschitz  hypersurfaces.

\proclaim{Lemma 4.1}    Let $M$ be a Lipschitz
Levi-flat hypersurface in   $\Bbb CP^n$, $n\ge 2$.  Then the
curvature form $\tilde\Theta^N$ associated
with the complex line bundle for $M$
  is  a well-defined 2-form with   $L^\infty$ coefficients in
$M$ and    is
$d$-exact. In fact, we have  $\tilde\Theta^N= \sqrt{-1}d\beta$ for some
  form $\beta$ on a
tubular neighborhood  $U (M)$
of
$M$ in
$\Bbb CP^n$. Furthermore, we can choose $\beta$ to be
$C^{1-\gamma}$-smooth for any small
$\gamma>0$.

\endproclaim
\demo{Proof}
  Let $Q$ be a point on $M$ and $\Sigma_Q$ be the holomorphic leaf of
$M$ passing through $Q$ with
  $\dim_{\Bbb C} [\Sigma_Q ] = n-1$. There is
  a holomorphic coordinate system $ (z_1, z_2, ..., z_n)  $ of $\Bbb
CP^n$ near $Q$, such that $
  (z_1, ..., z_{n-1}) $ is a local coordinate system of
$\Sigma_Q$ near $Q$. Applying the Gram-Schmidt process to the
local holomorphic frame $\{ \frac{\partial  }{\partial z_1},
\cdots,  \frac{\partial  }{\partial z_n}   \}$ near $Q$, we obtain
a special unitary basis
   $ \tilde e_1,\cdots,\tilde e_n$  such that  $\tilde e_l\in
T^{1,0}(M)$ for $l=1,\cdots,n-1$ and
$\tilde e_n |_P$  is orthogonal to $T_P^{(1, 0)}(\Sigma_Q)$ for
all $P \in \Sigma_Q$, with respect to the Fubini-Study metric. If
$M$ is $C^1$, then $\tilde e_n = \lambda
(\partial \rho)_\#$    for some $\lambda$ with $|\lambda|=\sqrt 2$.
Notice that $\lambda $ is
not necessarily a real valued function in $P
\in \Sigma_Q$. Let
  $\tilde\theta_{n,\bar l}$ be  the connection 1-forms with respect to
a unitary basis $\tilde
e_1,\cdots,\tilde e_n$  with $\tilde e_j \in T^{1,0}(M)$ for
$j=1,\cdots,n-1$. It is well-known that the curvature form
$\tilde\Theta^N$ of the quotient line bundle
$T^{(1,0)}(CP^n) / T^{(1,0)}( \Sigma_Q)$ is independent of the
choice of local frame $\{\tilde e_1,\cdots,\tilde e_n \}$.
Furthermore, its curvature form $\tilde\Theta^N$ is a closed form,
by the Chern-Weil theory. We remark that the Chern classes are well
defined for any continuous
complex vector bundle (see \cite{Mi}).

To see that $\tilde\Theta^N$ has $L^\infty$ coefficients, we use
the generalized Gauss-Codazzi equations (the Cartan-Chern
structure formula, see  (A.14)-(A.17) in \cite{CSW} and the
notation therein).
  Using Lemma 2.1, each
$\tilde \theta_{n,\bar l}$ has  bounded measurable coefficients.  Let
$\tilde \Theta$ denote the
curvature tensor for $\Bbb CP^n$  which is   an $n\times n$ matrix
and $\Theta_{n,\bar n}$ be its $(n,\bar n)$
component. We have
  $$\tilde \Theta^N=\tilde \Theta_{n,\bar n}-\sum_{l=1}^{n-1}\tilde
\theta_{n,\bar l}\wedge \tilde \theta_{l,\bar
n}, \tag 4.5$$ where $ \theta_{l,\bar n}  $  is given by
$$
\tilde \theta_{j,\bar n}( \xi) = - \langle \nabla_\xi (\tilde e_n
), \bar{\tilde e}_j  \rangle = \langle \nabla_\xi (\tilde e_j ),
\bar{\tilde e}_n  \rangle,\qquad \xi\in T(\Bbb CP^n),
$$
and $\theta_{j,\bar n} = - \bar{\theta}_{n,\bar j}  $  (see
(A.17)-(A.18) of [CSW]).
This gives that $\tilde \Theta^N$ has   bounded coefficients on  $M$.

Because $M$ has real codimension $1$ in $\Bbb CP^n$ and $M$ is
locally the  graph of some Lipschitz function, using a partition of
unity,   $M $ admits a nowhere vanishing
continuous global cross-section $\{\zeta\}$ in the quotient line
bundle $\Cal L = T^{(1, 0)}(\Bbb
CP^n) / T^{(1, 0)}(M)$. The quotient line
bundle $\Cal L$ is topologically trivial on $M$, just as in the
smooth case (see \cite{Si1}).

  This line bundle $\Cal L$ can be extended trivially to a small
neighborhood $U(M)$
of $M$.  Let
$\widehat M_s=\rho^{-1}(s)$. Then $\widehat M_s$ gives rise to a
family of Lipschitz hypersurfaces
for all
$|s|<\epsilon$, $\epsilon>0$ small,  with
$\widehat M_o=M$.   Using the mollifier smoothing
technique (cf
\cite{Ka}), one can obtain a family of smooth hypersurfaces $\tilde
M_s$ such that each $\tilde M_s$ is a
smooth real hypersurface when $s>0$ and  $\tilde M_0=M$. Let
$U_{\epsilon_0}=\cup_{|s|<\epsilon_0}\tilde M_s$
for some small
$\epsilon_0>0$.  Then
$U_{\epsilon_0}$ is an open neighborhood of $M$. The complex line
bundle $\tilde \Cal L$ on
$U_{\epsilon_0}$ induced by
  $T^{(1,0)}(\Bbb CP^n)/T^{(1,0)}(\tilde M_s)$ is
topologically trivial on $U_{\epsilon_0}$ since $\{ \zeta\}$ is a
nowhere vanishing continuous  cross
section.  Also    $\tilde \Cal L\mid _M=\Cal L$.    Thus the Chern
curvature form $\tilde \Theta^N$ is
$d$-exact  in $U_{\epsilon_0}$. Using (4.5) again, we see that
$\tilde \Theta^N$ has $L^\infty$ coefficients
in $U_{\epsilon_0}=U(M)$.

  Since $\tilde\Theta^N$ is $d$-exact on $U(M)$ and $L^\infty\subset
C^{-\gamma }$ on $U(M)$ for any $\gamma
>0$,  we can use
  the de Rham-Hodge decomposition theorem and interior regularity of the
$d$-operator on $U(M)$ to find some  $\beta $,  which is
$C^{1-\gamma}$-smooth for arbitrarily small $\gamma>0$. \qed
\enddemo

\proclaim{Proposition 4.2} Let $M$ be a compact Lipschitz Levi-flat
hypersurface in $\Bbb CP^n$,   $n\ge 3$. Let
$\beta_b^{0,1} $ be the projection of the Chern connection form
$\beta$ to $T^{0,1}(M)$, where $\beta\in
C^{1-\gamma}(M)$ is given by Lemma 4.1.   Then there exists an
$\epsilon'   > 0$ and a function $u\in C^{
\epsilon' }(M)$  such that
      $$\db_b u=\beta_b^{0,1}
\quad
\text{in } M .$$

\endproclaim
\demo{Proof} Let $\Bbb CP^n\setminus M=\Omega^+\cup\Omega^-$. Then
$\Omega^+$ and $\Omega^-$ are
pseudoconvex domains with
Lipschitz  Levi-flat boundary.
   From Lemma 4.1,     $     \beta_b^{(0,
1)} $ has   $C^{1-\gamma}$ coefficients where $0 <\gamma<1$ on $M$.

Since $M$ is Lipschitz, using the trace theorem (see [Gr]),  we can
extend $\beta_b^{0,1}$
to an $(0, 1)$-form  $\tilde \beta^{0,1}$on the whole $\Bbb CP^n$
such that $\tilde \beta^{0,1} \in W^{1 - \gamma + \frac 12}( \Bbb
CP^n  )$. Let  $0<\varepsilon =     \frac 12 -
\gamma  <\frac 12$. Then $f = \db \tilde \beta^{0,1}\in W^{s}(\Bbb
CP^n)$. We set $ f_\pm = f
|_{\Omega^\pm}$. We may choose our extension such that   that $f_\pm
\in L^2_{(0, 2)}(
\delta^{-t}, \Omega) $ for $t = 2\varepsilon$ since $\epsilon<\frac
12$.   Applying the proof
of Proposition 3.2 and Theorem  3.3, we observe that

$$
\hat \beta^{0,1}_\pm = \tilde \beta^{0,1} + \star_{(t)}
{\bar\partial N_t{\star_{(-t)} f}}_\pm
$$
is a $\db$-closed extension of $\beta_b^{0,1}$ to $\Omega^\pm$.
Thus, $\beta_b^{0,1}$ has a $\db$-closed extension $\hat
\beta^{0,1}$ on the whole $\Bbb CP^n$, with $\hat \beta^{0,1} \in
W^\varepsilon_{(0, 1)}(\Bbb CP^n)$. Since the cohomology group
$H^{(0, 1)}( \Bbb CP^n)$ vanishes. We can find   $u \in W^{1 +
\varepsilon}(\Bbb CP^n)$ with
$$
\db u = \hat \beta^{0,1}.
$$
Using the trace theorem   again, we conclude that there is a
$u\in W^{ \frac 12 + \epsilon}(M)$  such that

      $$\db_b u=  \beta_b^{0,1}\quad \text{on } M. \tag 4.6$$

      Using the local parametrization used in Lemma 2.1 with $V =
\cup_{|t| < \mu}
\Sigma_t \subset M$, the
equation
$\db_b $ is equal to $\db_{z'}$ on each leaf
$\Sigma_t$, which is elliptic. From Lemma 4.1,
$\beta_b^{0,1}$ is $C^{\alpha}$ on $M$.
    From (4.6),
and the classic Schauder theorem (cf. [GT]) for elliptic equations on
$\Sigma_t$,
we get that  $u$ is
$C^{1,\alpha}$-smooth on each leaf.  Furthermore, we have (see e.g.
\cite{ShW}) that
there exists a constant $C $   independent of
$t$    such that
$$|u|_{C^{1+\alpha}(\Sigma_t)}\le C
(|\beta_b^{0,1}|_{C^{\alpha}(\Sigma_t)}+\|u\|_{L^2 (\Sigma_t)}),\tag
4.7$$ where $C$ depends on the neighborhood $V$ of $Q$ and the
parametrization   $\Psi$, but   is
independent of
$t$ since (4.6) is uniformly elliptic on $\Sigma_t \subset V$
independent of $t$.

   From the Sobolev
trace theorem,  the function $u\in W^{\frac 12+ \epsilon}(M) $ has
$L^2$-trace on each leaf. Therefore,
      there exists $C_2>0$  independent of $t$ such that
$$\|u\|_{L^2 (\Sigma_t)}\le C_2 \|u\|_{W^{\frac 12+ \epsilon} (M)}.
\tag 4.8$$
Combining (4.7) and (4.8), we get
$$   |u|_{L^\infty(V)}\le \underset{|t| < \mu}\to\sup
|u|_{C^{1+\alpha}(\Sigma_t)}\le C_3.\tag 4.9 $$
Thus we have  already proved that  $u$ is bounded.

   It remains to prove that
$u$ is H\"older continuous in the transversal
$t$ direction.  We can prove this by  applying  a
modified one-dimensional Sobolev embedding theorem. This can be done
by  taking the finite difference
of  the equation (4.6) with respect to  the Besov norms.
     The proof is exactly the same as
       before and we refer the reader to the proofs of Lemmas
5.2-5.3 in
\cite{CSW}. Thus we conclude that
     $u\in C^{ \epsilon' }(M)$ for some sufficiently small $\epsilon'
< \epsilon$.\qed
\enddemo

\demo{Proof of the  theorem} Using Lemma 4.1 and (4.5), we have that
the curvature form
$\sqrt{-1}\tilde\Theta_b^N$ is positive definite on each holomorphic
leaf of  the   Levi-flat hypersurface $M$
(see Proposition A.2 in the  Appendix in \cite{CSW}).
Let $h=2\text{Im}u$, where $u$ is the function obtained in
Proposition 4.2.  We have
$$\sqrt{-1}\partial_b\db_b h=\sqrt{-1} \tilde \Theta^N_b>0
\quad\text{on } {T^{(1, 0)}(M)\oplus T^{(0,
1)}(M) }.  \tag 4.10$$ Since $h$ is continuous on the compact
hypersurface  $M$, it attains
its maximum at some point $p$  in $M$. Since $p$ lies  in the interior of some
leaf, one   obtains a contradiction from (4.10) and  the    Maximum
Principle. This completes the proof of the
theorem.
\qed
\enddemo

\heading {\bf 5. The case for   $\Bbb CP^2$ }\endheading

To prove the nonexistence of Levi-flat hypersurfaces in $\Bbb CP^2$,
we can study the
$\db$-Cauchy problem for the top degree forms.
There are major differences for compatibility conditions for
$\db$-closed extensions of $(0, q)$-forms
when  $q <n-1$ and   $ q=n-1 $.
In general,  the space of  harmonic $(p,n-1)$-forms on a
pseudoconcave domain in $\Bbb CP^n$
is infinite dimensional (see Theorem 3.1 in H\"ormander \cite{H\"o2}).

For $q=n-1$, there is an additional compatibility
condition  for the $\db$-closed  extension of
$(p,n-1)$-forms.

\proclaim{Proposition 5.1} Let $\Omega \subset\subset \Bbb CP^n$ be a
pseudoconvex domain with $C^2$
boundary, $n\ge 2$, and let $\Omega^+=\Bbb CP^n\setminus \overline\Omega$.
   For any $\db$-closed $f\in W^1_{(p,n-1)}( {\Omega}^+)$, where
$0\le p\le n$, the following conditions are equivalent:
\roster
\item    There exists
$F\in L^2_{(p,n-1)}(\Bbb CP^n)$ such that
$F|_{\Omega^+}=f$ and  $\db F=0$ in $\Bbb CP^n$ in the
distribution sense.

   \item    The restriction of  $f$ to $b\Omega$  satisfies the
compatibility condition
$$
\int_{b\Omega}    f\wedge \phi     =0,  \quad\phi\in L^2_{(n-p,0)} (\Omega)
\cap\text{Ker}(\db).  $$
\item  Any $W^1 $ extension $\tilde f\in W^1_{(p,n-1)}(\Bbb CP^n) $
of  $f$ satisfies the
compatibility condition
$$
\int_{ \Omega}   \db\tilde f\wedge \phi     =0,  \quad\phi\in
L^2_{(n-p,0)} (\Omega)
\cap\text{Ker}(\db).  $$
   When $p\neq n-1$, the above conditions are equivalent to
\item There exists $u\in W^1_{(p,n-2)}(\Omega^+)$ satisfying
   $\db u=f$ in $\Omega^+. $
\endroster

\endproclaim
   We remark that  any $f$ in $W^1(\Omega^+)$ has a trace in $W^{\frac
12}(b\Omega^+)$ and any
holomorphic
$(n-p,0)$-form with $L^2(\Omega)$ coefficients  has a well-defined
trace in $W^{-\frac 12}(b\Omega)$
(see e.g. \cite{LM}). Thus the pairing between $f$  and
$\phi$   in (2) is well-defined.

\demo{Proof} We first show  that (1) implies (2).

We assume that  there exists a $\db$-closed extension $F$ of a
$\db$-closed form $f$.   For
    $\phi\in C^1_{(n-p,0)} (\overline\Omega)
\cap\text{Ker}(\db)$,
by Stokes' theorem,   we have
$$ 0=\int_{\Omega}  \db F\wedge \phi=\int_{\Omega}  \db (F\wedge \phi)=
\int_{b\Omega}     f\wedge \phi  =0.$$
If $\phi$ is only in $L^2$, we use an approximating sequence
$\phi_\nu\in  C^1_{(n-p,0)}
(\Omega)$ such that $\phi_\nu\to \phi$ in $ L^2_{(n-p,0)} (\Omega) $ and
$\db \phi_\nu\to 0$ in $ L^2_{(n-p,1)} (\Omega) $ by the Friedrichs'
Lemma (cf. [CS]).
We have
$$\aligned  0&=\lim_{\nu\to \infty}(\int_{\Omega}  \db F\wedge
\phi_\nu+(-1)^{p+n-1}\int_{\Omega}
F\wedge\db
\phi_\nu)\\&=\lim_{\nu\to \infty}
\int_{b\Omega}     f\wedge \phi_\nu=
\int_{b\Omega}     f\wedge \phi   .\endaligned$$

To see that (2) implies (3), we observe
$$\int_{ \Omega}   \db\tilde f\wedge \phi =\int_{ b\Omega}    f\wedge
\phi     =0    .$$

To show that (3) implies (1), we set
$$  T\tilde f = - {\star\bar\partial N_{(n-p,0)} {(\star  \db \tilde f)}}
\quad \quad
\text{    on } \Omega.
$$
  From the proofs of  Theorem 1.2 or Corollary 1.6, we have
$\db T\tilde f= \db \tilde f$
in $\Bbb CP^n$ if we extend $T\tilde f$ to be zero outside $\Omega$.
Define $F$ the same as in (2.1).
   Then $F\in L^2_{(p,n-1)}(\Bbb CP^n)$ and $F$ is a $\db$-closed
extension of $f$. This proves
that conditions (1), (2) and (3) are equivalent.

When $p\neq n-1$, the harmonic $(p,n-1)$-forms  $\Cal H_{(p,n-1)}(\Bbb
CP^n)=\{0\}$. Thus if (1)
holds, then there exists  $u\in W^1_{(p,n-2)}(\Bbb CP^n)$ satisfying
    $\db u=f$ in $\Bbb CP^n$. Restricting $u$ to $\Omega^+$, we have proved (4).
Conversely, if
$f$ is
$\db$-exact for some
$u\in W^1_{(p,n-2)}(\Omega^+)$, we can extend $u$ to be a
$(p,n-2)$-form in $W^1_{(p,n-2)}(\Bbb
CP^n)$. Then the $(p,n-1)$-form
$F=\db u$ is a
$\db$-closed extension of $f$ with $L^2$ coefficients. Thus (1) and
(4) are equivalent.  The
    proposition is proved.
\qed
\enddemo

Proposition 5.1 also holds for any $\Omega$ with $C^{1, 1}$ Levi-flat boundary.
If one   can show  that
any $\db$-closed form  on $\Omega^+$ with $W^1(\Omega^+)$ coefficients extends
to be $\db$-closed
in $\Bbb CP^2$, i.e., any of the equivalent conditions in Proposition
5.1 holds on a domain with
$C^{1, 1}$ Levi-flat boundary,  then one can show  the nonexistence of
$C^{1, 1}$ Levi-flat hypersurfaces in
$\Bbb CP^2$  using arguments similar to the proof of the main theorem
in Section 4.
Notice that in this case, the domain $\Omega$ is both pseudoconvex
and pseudoconcave.
But to prove the nonexistence of  Levi-flat hypersurfaces in $\Bbb CP^2$,  we
need the following
$W^1$ regularity for the $\db$-equation.
\proclaim{Conjecture 1 ($W^1$ regularity for $\db$)}   Let $\Omega
\subset\subset\Bbb CP^2$ be a
      Lipschitz    domain with    Levi-flat
boundary.  For any    $f\in
C^\infty_{(0,1 )}(\overline{\Omega})$ with $\db f=0$,  there exists
       $u\in W^1 (\Omega)$ such
that $\db u=f$.
\endproclaim

    Conjecture 1 will yield the nonexistence of Lipschitz  Levi-flat
hypersurfaces in
$\Bbb CP^2$.   When $b\Omega$ is
$C^4$ and Levi-flat,    this is proved by Siu (see \cite{Si2}) with
$u\in W^3(\Omega)$. It seems that one only
needs  the boundary to be $C^2$ to have a solution $u\in W^1
(\Omega)$. Thus we can reduce the smoothness
assumption used in \cite{Si2} on $\Omega$, but the $W^1$ regularity
of the solution for the
$\db$-equation cannot be removed.

   The following Liouville type result   stated in Proposition 4.5
in \cite{CSW} remains open.
\proclaim{Conjecture 2 (Liouville's Theorem)} Let
$\Omega^+\subset\subset\Bbb CP^n$ be a pseudoconcave domain with
$C^2$-smooth boundary (or
Lipschitz) $b\Omega^+ $, $n\ge 2$.   Then
$L^2_{(p,0)}( \Omega^+)\cap \text{Ker}(\db)=\{0\}$  for every $1\le  p\le n$
and  $L^2_{(0,0)}( \Omega^+)\cap \text{Ker}(\db)=\Bbb C$.
\endproclaim
This conjecture also implies  the nonexistence of Levi-flat
hypersurfaces in $\Bbb CP^n$
for $n\ge 2$.  From Corollary 3.4, the set $W^{1+\epsilon}_{(p,0)}(
\Omega^+)\cap \text{Ker}(\db)
$ is either zero or constants  for Lipschitz pseudoconcave domains.
When the boundary is $C^2$, this is also
true for
$\epsilon=0$.   Thus it suffices to show that
$W^{1}_{(p,0)}( \Omega^+)\cap \text{Ker}(\db)
$ is dense in $L^2_{(p,0)}( \Omega^+)\cap \text{Ker}(\db)$ for the $C^2$ case.
There is  still a gap in the  the required uniform estimates
(4.18) for Proposition 4.5  in
\cite {CSW}. We remark that Conjecture 2 is much stronger than the
nonexistence of Levi-flat
hypersurfaces, since there are many pseudoconcave domains in $\Bbb CP^n$.
\bigskip
\noindent
{\bf Acknowledgment}. We would like to
 thank   Sophia Vassiliadou for pointing out some errors
in
Section 2 in the original manuscript.


\Refs
\nofrills{References}
\widestnumber\key{APS12}

\ref \key BF \by Barrett, D. E. and  Fornaess, J. E. \paper On the smoothness
of Levi-foliations\jour  Publ. Mat. \vol 32  \yr 1988\pages  171-177\endref

\ref \key BC \by Berndtsson, B and Charpentier, P. \paper
A Sobolev mapping property of the Bergman kernel \jour  Math. Zeitschrift\vol
235
\yr 2000 \pages 1-10\endref

\ref\key CS\by Cao, J. and Shaw, M.-C. \paper A new proof of the
Takeuchi theorem
\jour preprint\endref

\ref\key  CSW  \by   Cao, J.,  Shaw,  M.-C. and   Wang L.
\paper Estimates for the $\db$-Neumann problem     and nonexistence of
         Levi-flat hypersurfaces  in $\Bbb CP^n$
\jour Math. Zeit \vol 248
\yr 2004 \pages 183-221  \moreref \paper \rom{Erratum}\pages 223-225
\endref






\ref \key CS \by Chen, S.-C. and Shaw, M.-C. \book Partial Differential
Equations
in Several Complex Variables
\publ American Math. Society-International Press, Studies in Advanced
Mathematics, Volume 19\publaddr Providence, R.I.
\yr 2001\endref

\ref\key De\by Demailly, J.-P. \paper Mesures de Monge-Amp\`ere et
mesures pluriharmoniques\jour
Math. Zeit.\vol 194\pages 519-564\yr 1987\endref

\ref\key DF\by Diederich, K.  and Fornaess, J. E.\pages
129--141
\paper Pseudoconvex domains: Bounded strictly
plurisubharmonic exhaustion functions
\yr 1977\vol 39\jour Invent. Math.,
\endref

\ref \key EG\by Evans, L. E.  and  Gariepy, R. F. \book {\rm Measure theory and
fine properties of functions}\publ
CRC press\publaddr Boca Raton \yr 1992\endref



\ref\key Fe\by Federer, H.\paper Curvature measures \jour Trans.
Amer. Math. Society
\vol 93\yr 1959\pages 418-491
\endref

\ref\key FK \by Folland, G. B.  and Kohn, J. J.
\book {\rm The Neumann Problem for the Cauchy-Riemann
Complex, Ann. Math. Studies 75}
\yr 1972\publ Princeton University Press, Princeton, N.J.
\endref


\ref \key Gr \by Grisvard, P.
\book {\rm Elliptic Problems in Nonsmooth Domains}
\publ Pitman\publaddr Boston\yr 1985
\endref


\ref\key HI \by Henkin, G. M. and Iordan, A. \paper Regularity of $\db$ on
pseudoconcave compacts and applications \jour Asian J. Math. \vol
4\yr 2000\pages 855-884 \paperinfo  Erratum: Asian J. Math., vol
{\bf 7}, (2003) No. 1,  pp. 147-148) \endref










\ref  \key H\"o1 \by   H\"ormander, L. \paper $L^2$ estimates
and existence theorems for the $\bar\partial$
operator \jour Acta Math. \vol  113 \yr 1965\pages
89-152\endref

\ref  \key H\"o2  \by   H\"ormander, L. \paper  The null space of the
$\overline\partial$-Neumann
operator  \jour Ann. Inst. Fourier
(Grenoble)   \vol  54\yr 2004\pages 1305-1369\endref




\ref\key Io \by Iordan, A.\paper On the non-existence of smooth Levi-flat
hypersurfaces in $\Bbb CP^n$\paperinfo will appear in the
``Proceedings of the Memorial Conference of Kiyoshi Oka's Centenial
Birthday on Complex Analysis in Several Variables", Kyoto, Nara 2001
\endref

\ref \key Ka\by Karcher, H.\paper Riemannian center of mass and
mollifier smoothing\jour Comm. Pure and
Appl. Math. \vol 30\yr 1977\pages 509-541\endref


\ref\key KeR\by Kerzman, N. and Rosay, J.-P.\paper Fonctions
Plurisousharmoniques d'exhaustion born\'ees et domaines
taut\jour Math. Ann. \vol 257\pages 171-184\yr 1981\endref

\ref\key KR   \by Kohn, J. J., and Rossi, H.\pages 451-472
\paper On the extension of holomorphic functions from the
boundary of a complex manifold
\yr 1965\vol 81\jour Ann. Math.,
\endref


\ref\key LM  \by Lions, J.-L., and Magenes, E.
\book {\rm Non-Homogeneous Boundary Value Problems and
Applications, Volume I}
\publ Springer-Verlag, New York \yr 1972
\endref



\ref \key LN  \by  Lins Neto, A. \paper A note on projective Levi flats and
minimal sets of algebraic foliations
\jour Ann. Inst. Fourier\vol 49\yr 1999\pages 1369-1385
\endref

\ref\key Mi\by Milnor, J. W.\book {\rm Characteristic classes}
\publ Princeton University Press, Princeton, N. J.  \yr 1974
\endref







    \ref \key  NW \by Ni, L. and   Wolfson, J.\paper The Lefschetz
theorem for CR submanifolds and
the nonexistence of real analytic Levi
      flat submanifolds\jour  Comm. Anal. Geom.\vol 11 \yr 2003\pages 553-564
\endref


















\ref\key OS \by Ohsawa, T. and  Sibony, N.  \paper
Bounded P.S.H Functions and Pseudoconvexity in K\"ahler Manifolds \jour
Nagoya Math.
J.\vol 149\yr 1998\pages 1-8\endref



\ref\key ShW \by   Shaw,   M.-C. and  Wang, L \paper H\"older and $L^p$
estimates   for $\square_b$ on CR manifolds of arbitrary codimension
\jour Math. Ann. \vol 331\yr 2004\pages 297-343
\endref








\ref\key Si1 \by Siu, Y.-T. \paper Nonexistence of smooth Levi-flat
hypersurfaces in complex projective spaces of
dimension
$\ge 3$\jour Ann. Math. \vol 151\yr 2000\pages 1217-1243
\endref


\ref\key Si2 \by Siu, Y.-T.\paper $\db$-regularity for weakly pseudoconvex
domains in hermitian symmetric spaces with
respect to invariant metrics  \jour   Ann. Math.\vol 156\pages  595-621 \yr
2002
\endref

\ref\key St \by Stein, E. M.
\book {\rm Singular Integrals and Differentiability
Properties of Functions, Math. Series 30}
\publ  Princeton University Press, Princeton, New Jersey\yr
1970
\endref




\ref \key Ta \by Takeuchi A. \paper Domaines pseudoconvexes infinis et la
m\'etrique riemannienne dans un espace
projectif\jour J. Math. Soc. Japan\vol 16\yr 1964\pages 159-181\endref
\smallskip


\endRefs
\enddocument

\end